\documentclass[twocolumn]{autart}    

\usepackage{mathtools}
\mathtoolsset{showonlyrefs,showmanualtags}

\usepackage{datetime}
\ddmmyyyydate
\settimeformat{xxivtime}
\usepackage{amsmath}
\usepackage{amssymb}

\usepackage{color}
\usepackage{url}

\usepackage[abbrev]{bib_sty} 

\usepackage{paralist}

\usepackage{graphicx}

\usepackage{pgfplots} 
\pgfplotsset{compat=newest} 
\pgfplotsset{plot coordinates/math parser=false} 
\newlength\figureheight 
\newlength\figurewidth 

\newcommand{\Fig}{Figure}
\newcommand{\Sec}{Section}

\newcommand{\Hinf}{\ensuremath{\mathcal{H}_\infty}}
\newcommand{\lambdamax}{\ensuremath{\lambda^{\mathrm{max}}}}

\newcommand{\iter}{j}
\newcommand{\System}{P}
\newcommand{\tijd}{k}
\newcommand{\tijdreverse}{\mathcal{T}_N}
\newcommand{\lambdashift}{\tilde \lambda}

\newtheorem{lemma}[thm]{Lemma}

\newcommand{\Colorone}{blue}
\newcommand{\Colortwo}{red}
\newcommand{\Colorthree}{black}

\usepackage{eso-pic}
\AddToShipoutPictureBG*{%
\AtPageUpperLeft{%
\setlength\unitlength{1in}%
\hspace*{\dimexpr0.5\paperwidth\relax}
\makebox(8,-.75)[c]{
\begin{tabular}{c c}
Tom Oomen and Cristian R. Rojas,\\
Reset-Free Data-Driven Gain Estimation: Power Iteration using Reversed-Circulant Matrices, \\
to appear in
{\em Automatica},
uploaded to ArXiv \today \\
\end{tabular}}}}

\begin{document}

\begin{frontmatter}

\title{Reset-Free Data-Driven Gain Estimation:\\ Power Iteration using Reversed-Circulant Matrices}

\author[Oomen]{Tom Oomen}\ead{t.a.e.oomen@tue.nl},    
\author[Rojas]{Cristian R. Rojas}\ead{crro@kth.se},               

\address[Oomen]{Department of Mechanical Engineering, Eindhoven University of Technology, Eindhoven, The Netherlands and Delft Center for Systems and Control, Delft University of Technology, Delft, The Netherlands.}  
\address[Rojas]{Division of Decision and Control Systems, KTH Royal Institute of Technology, 100 44 Stockholm, Sweden.}             

\begin{keyword}                           
Identification and control methods;
Data-driven control;
Data-driven robust control;
Identification for control;
Input and excitation design.
\end{keyword}

\begin{abstract}                          
A direct data-driven iterative algorithm is developed to accurately estimate the $\Hinf$ norm of a linear time-invariant system from continuous operation, i.e., without resetting the system. The main technical step involves a reversed-circulant matrix that can be evaluated in a model-free setting  by performing experiments on the real system.
\end{abstract}

\end{frontmatter}

\section{Introduction}

Direct data-driven approaches to estimate system properties, such as system norms, enable the extraction of information while avoiding the need to model the entire system. In \cite[\Sec\ 12.2]{Hjalmarsson2005}, an iterative approach to estimate the $\Hinf$ norm directly from data is developed, which is extended and analyzed in \cite{RojasOomHjaWah2012}, \cite{TuBocRec2018}, \cite{RojasMul2019}, \cite{KochMonAll2021}. In sharp contrast, model-based approaches first estimate a model, either a parametric model with subsequent norm calculation \cite{BruinsmaSte1990}, or a nonparametric model \cite{GeerardynOom2017}. These model-based approaches require substantial user intervention compared to data-driven approaches. Furthermore, these data-driven approaches can be viewed as iterative experiment design, see \cite[\Sec\ 4.2]{RojasOomHjaWah2012}. 

A key assumption that underlies the developments in these direct data-driven norm estimation approaches is that the system needs to be reset after each iteration, which tacitly leads to substantial theoretical and practical drawbacks. The main reason to assume a full reset is that it allows to write the system as a linear system of equations, enabling a direct application of a power iteration. If such a reset is not performed, then past iterations still have effect during the current iteration, which is further complicated by a nonlinear normalization \cite{RojasOomHjaWah2012} at each iteration. A reset directly leads to a theoretical drawback, since the convergence rate of the estimated norm in these reset-based approaches is $O(N^{-2})$  \cite{TuBocRec2019}, where $N$ denotes the experiment length of a single iteration. A frequency-domain sampling approach would lead to faster convergence results by eliminating transient effects, see \cite{GeerardynOom2017}, but still requires resets and does not benefit from the earlier mentioned optimal experiment design compared to iterative schemes. As major practical drawback, the application of the required resets requires additional experimental time and effort, which is application dependent but typically intensive.

Although many approaches have been developed that enable data-driven estimation of system properties, these approaches typically require the system to be reset every iteration. The aim of this paper is to develop an approach that has improved accuracy and facilitates practical implementation by avoiding resets.

The main contribution of this paper is a reset-free iterative data-driven algorithm. The main technical result is a periodic response matrix that is transformed into a new reversed-circulant matrix. Specific focus is on a $\mathcal{H}_\infty$-norm estimation procedure, similar to \cite[\Sec\ 12.2]{Hjalmarsson2005}, \cite{RojasOomHjaWah2012}, \cite{TuBocRec2018}, \cite{RojasMul2019}, \cite{KochMonAll2021} yet reset-free. In particular, by analyzing the eigenvalue decomposition of the reversed-circulant matrix, an iterative power iteration method is devised that estimates the $\mathcal{H}_\infty$-norm using real-valued inputs. 

The main direct benefit of the proposed method is a data-driven algorithm that provides more accurate estimates of the $\mathcal{H}_\infty$-norm compared to reset-based approaches, which is confirmed in a simulation example. In addition, there is a substantial benefit in practical applications, including the one reported in \cite{OomenMaaRojHja2014}, where the system does not need to be reset through a cumbersome re-initialisation and homing procedure, reducing experimental time and user effort. The broader significance of the paper lies in advocating the possibility that reset-free, i.e., continuous operation, may have theoretical and practical advantages for general data-driven estimation schemes, cf.~\cite{KochMonAll2021} for a recent overview, which may unnecessarily restrict to reset-based approaches that are unnecessarily theoretically conservative and experimentally impractical. Note that these benefits are recognised in certain specific control approaches, e.g. in \cite{RozarioOom2023} and references therein, where also continuous operation is advocated. It is thus expected that continuous operation may have benefits for broader classes of data-driven algorithms, e.g., \cite{Hjalmarsson2005}, \cite{BazanellaCamEck2012}.

\section{Problem definition}

\subsection{Iterative experiments}
Consider the discrete-time single-input single-output linear time-invariant system $\System(z)$ given by the state-space realization
\begin{align}\label{eq:ss}
x(\tijd+1) &= A x(\tijd) + B u(\tijd) \\
y(\tijd) &= C x(\tijd) + D u(\tijd),
\end{align}
with matrices of appropriate dimensions and $\lambdamax(A)<1$. Data are collected in batches of length $N \in \mathbb{Z}_{>0}$, where the $\iter^{\mathrm{th}}$ batch or iteration is denoted by
\begin{equation}
y_\iter 
=
\begin{bmatrix}
y(N\iter) & y (N\iter + 1) & \ldots & y (N\iter + N - 1)
\end{bmatrix}^T.
\end{equation}
Consequently, by processing $N$ data batch-wise,
\begin{align}\label{eq:liftedsystem}
x(N\iter + N) &= F x(N\iter) + G u_\iter\\
y_\iter &= H x(N \iter) + J u_\iter
\end{align}
is obtained. Here,
\begin{align}\label{eq:FGHJ}
\left[
\begin{array}{c|c}
F & G \\\hline H & J
\end{array}
\right]
= 
 \left[
\begin{array}{c|ccccc}
A^N & A^{N-1} B &  \ldots & AB & B \\\hline 
C & D  & \cdots & 0 & 0 \\
CA & CB  & \cdots & 0 & 0 \\
\vdots & \vdots & \ddots & \ddots & \vdots \\
CA^{N-1} & CA^{N-2} B  & \cdots  & CB & D
\end{array}
\right].
\normalsize
\end{align}

\subsection{Traditional data-driven power iterations via resetting}\label{sec:PIM}
In this section, the main idea behind traditional data-driven norm estimation is outlined, including \cite{RojasOomHjaWah2012}. The system \eqref{eq:liftedsystem} is reset after each experiment, i.e., 
\begin{equation}\label{eq:reset}
x(0) = x(N) = x(2N) = \cdots = x(N\iter) = 0,
\end{equation}
hence \eqref{eq:liftedsystem} collapses to
\begin{equation}
y_\iter = J u_\iter.
\end{equation}
Next, the induced $2$-norm is computed as 
\begin{equation}\label{eq:Ri2_fin_time}
  \|J\|_{i2} := \sup_{u \neq 0}\tfrac{\|y\|_{2}}{\|u\|_2}  = \sup_{u_\iter \neq 0} \sqrt{\tfrac{u_\iter^TJ^TJu_\iter}{u_\iter^Tu_\iter}}
  =  \sqrt{\lambdamax(J^TJ)}.
\end{equation}
Note that
\begin{equation}\label{eq:transposeistimereversals}
  J^T = \tijdreverse J \tijdreverse,
\end{equation}
 where $\tijdreverse \in \mathbb{R}^{N \times N}$ is the involutory permutation matrix 
\begin{equation}
\tijdreverse = \begin{bmatrix}
0 & \cdots & 0 & 1\\
0 & \cdots & 1 & 0\\
\vdots &\ddots & \ddots & \vdots\\
1 & \cdots & 0 & 0
\end{bmatrix},
\end{equation}
which can be interpreted as a time-reversal operator. Next, \eqref{eq:transposeistimereversals} reveals that $\tijdreverse J$ is symmetric, hence
\begin{equation}\label{eq:Ri2_fin_time3}
  \|J\|_{i2} =  {\lambdamax(\tijdreverse J)}.
\end{equation}
As a result, $\|J\|_{i2}$ equals the largest eigenvalue $\lambdamax(\tijdreverse J)$. In addition, $\lambdamax(\tijdreverse J) \rightarrow \|P\|_\infty$ for $N \rightarrow \infty$, see \cite[Theorem 3]{RojasOomHjaWah2012} for a proof, where for SISO systems
\begin{equation}
\|P\|_\infty = \sup_{\omega}\left| P(e^{j\omega})\right|.
\end{equation}
The main idea of the data-driven approach is enabled by \eqref{eq:Ri2_fin_time3}. In particular, \eqref{eq:Ri2_fin_time3} enables that all knowledge of the system $J$ can directly be evaluated on the system instead of using model knowledge. In particular, the power method \cite[\Sec\ 9.1]{Wilkinson1965}, which is an iterative algorithm to determine the maximal eigenvalue of a matrix, can be directly applied to the real system $J$ by applying a certain input and measuring the output. Then, the output is time-reversed through the operation $\tijdreverse$ and re-applied as input. By initializing with a random $u_0$, the input  $\lim_{\iter \rightarrow \infty} u_\iter$ converges to an eigenvector corresponding to $\lambda_{\mathrm{max}} = \max |\lambda(J \tijdreverse)|$. As a result, the maximum eigenvalue directly follows from a Rayleigh quotient.

\subsection{Problem formulation}
Traditional iterative methods in \Sec~\ref{sec:PIM} are based on the reset in \eqref{eq:reset} and have several shortcomings. First, $\lim_{N \rightarrow \infty} \lambda_{\mathrm{max}} (J \tijdreverse) = \|\System\|_\infty$ is a limit result, with a relatively slow convergence rate of $O(N^{-2})$, where $N$ denotes the experiment length, see, e.g., \cite{TuBocRec2019} for a recent investigation of this and \cite{GeerardynOom2017,MullerRoj2020} for approaches that eliminate such transients effects in nonparametric model-based and model-free settings, respectively. Second, the reset requires substantial experimental effort to stop the system completely. The main problem addressed in this paper is to develop a direct data-driven $\Hinf$-norm estimation procedure for the system \eqref{eq:liftedsystem} that alleviates the shortcomings associated with resetting \eqref{eq:reset}. 

\section{Reset-free data-driven estimation approach}

The main idea is to devise a data-driven mechanism that excites the system \eqref{eq:liftedsystem} with a periodic input, i.e., 
\begin{equation}\label{eq:convinput}
u_\iter = u_{\iter+1} = \cdots = u_\infty.
\end{equation}
Assuming $\rho(A)<1$, hence $\rho(F)<1$, this leads to
\begin{equation}\label{eq:periodicresponse}
x_\infty = F x_\infty + G u_\infty
= (I-F)^{-1}Gu_\infty
\end{equation} 
and thus
\begin{equation}\label{eq:HIFGJ}
y_\infty = H x_\infty + J u_\infty = (H(I-F)^{-1}G+J)u_\infty.
\end{equation}
Hence, continuous reset-free operation essentially enables experimenting on the $N \times N$ periodic response matrix $H(I-F)^{-1}G+J$. The maximum gain of this matrix corresponds to the $\Hinf$ norm for $N \rightarrow \infty$, as is revealed next.

\begin{thm}\label{thm:circper}
Consider the matrix $H(I-F)^{-1}G+J$ defined by \eqref{eq:liftedsystem}, and let
\begin{equation}\label{eq:Fourier}
\mathcal{F} = \frac{1}{\sqrt{N}} \begin{bmatrix}
1 & 1 & 1 & \cdots & 1 \\
1 & W_N & W_N^2 & \cdots & W_N^{N-1}\\
1 & W_N^2 & W_N^4 & \cdots & W_N^{2(N-1)}\\
\vdots & \vdots & \vdots & \ddots & \vdots \\
1 & W_N^{N-1} & W_N^{2(N-1)} & \cdots & W_N^{(N-1)(N-1)} 
\end{bmatrix}
\end{equation}
be the $N \times N$ DFT matrix, with $W_N = e^{-\frac{2\pi j}{N}}$ a primitive $N^\mathrm{th}$ root of unity, i.e., satisfying $(W_N)^N = 1$. Then, 
$\mathcal{F}^{*}(H(I-F)^{-1}G+J)\mathcal{F} = \mathrm{diag} (G(e^{j\omega_m}))$.
\end{thm}
Note that the DFT matrix $\mathcal{F}$ is unitary, i.e., $\mathcal{F}^* \mathcal{F} = \mathcal{F} \mathcal{F}^* = I$. Next, circulant matrices are defined, since these matrices play a key role in the proof of Theorem~\ref{thm:circper} and the main procedure of this paper.
\begin{defn}
Let $N \geq 1$ and $a_0, a_1, \ldots a_{N-1} \in \mathbb{R}$. Then,
\begin{multline}\label{eq:circ}
\mathrm{circ}(a_0, a_1, \ldots a_{N-1}) 
= \\
 \begin{bmatrix}
a_0 & a_1 & a_2 & \cdots & a_{N-2} & a_{N-1} \\
a_{N-1} & a_0 & a_1 & \cdots & a_{N-3} & a_{N-2}\\  
a_{N-2} & a_{N-1} & a_0 & \ddots &  & a_{N-3}\\ 
\vdots & \vdots & \ddots & \ddots & \ddots & \vdots \\
a_2 & a_3 & a_4 & \ddots & a_0 & a_1 \\
a_1 & a_2 & a_3 & \cdots & a_{N-1} & a_0 
\end{bmatrix}
\end{multline}
is a \emph{circulant matrix}.
\end{defn}
\begin{lemma}\label{lemma:Jpervalues}
The matrix $H(I-F)^{-1}G+J$ in \eqref{eq:HIFGJ} with $F,G,H,J$ defined in \eqref{eq:FGHJ} is circulant, with 
\begin{equation}
a_k = 
\begin{cases}
D + C A^{N-1}(I-A^N)^{-1}B, & k = 0 \\
C A^{N-k-1}(I-A^N)^{-1}B, & k = 1, \ldots, N-1,
\end{cases}  \label{eq:a_k}
\end{equation}
eigenvalue decomposition
\begin{math}
\mathcal{F}^{-1}\mathrm{circ}(a_0, a_1, \ldots a_{N-1})\mathcal{F} = \mathrm{diag}(\lambda_0, \lambda_2, \ldots,\lambda_N-1)
\end{math}, and
\begin{equation}\label{eq:dft}
\lambda_m = \sum_{k=0}^{N-1} a_k e^{-\frac{2\pi jmk}{N}},\ m = 0, 1, \ldots, N-1.
\end{equation}
\end{lemma}
\begin{pf}[Lemma~\ref{lemma:Jpervalues}]
For the lower-triangular part, where $J$ is possibly nonzero, the $(p,q)^{\mathrm{th}}$ element ($p>q$) of $H(I-F)^{-1}G+J$ is
\begin{math}
CA^{p-1} (I-A^N)^{-1} A^{N-q}+ CA^{p-q-1}B 
= CA^{p-1} (I-A^N)^{-1} A^{N-q}B + CA^{p-q-1} (I-A^N)(I-A^N)^{-1}B 
= CA^{p-1} (I-A^N)^{-1} A^{N-q}B - CA^{p-q-1} A^N(I-A^N)^{-1}B 
+ CA^{p-q-1}(I-A^N)^{-1}B, \label{eq:result1}
\end{math}
where use is made of $(I-A^N)^{-1}A = A (I-A^N)^{-1}$, which can be verified by premultiplying and postmultiplying both sides by $I-A^N$. For the upper-triangular part, the $(p,q)^{\mathrm{th}}$ element ($q>p$) of $H(I-F)^{-1}G+J$ is given by
\begin{math}
CA^{p-1}(I-A^N)^{-1}A^{N-q}B 
= CA^{N+p-q-1}(I-A^N)^{-1}B.
\end{math} 
Comparing these entries to those in $\eqref{eq:circ}$ reveals that they indeed lead to a circulant matrix. Next, \eqref{eq:a_k} directly follows by noting that $a_k$ is at the $(1,i+1)$ location of \eqref{eq:circ} for $k>0$, while for $a_0$ an additional $D$ term is included. Finally, any circulant matrix satisfies the eigenvalue decomposition \eqref{eq:dft}, see
\cite[Fact 5.16.7]{Bernstein2005}. 
\end{pf}
\begin{rem}
In Lemma~\ref{lemma:Jpervalues}, circularity of $H(I-F)^{-1}G+J$ follows from the fact that \eqref{eq:FGHJ} is LTI. Lifting in \eqref{eq:liftedsystem} is often also performed for linear periodically time varying systems, in which case the matrix $\mathcal{F}^{*}(H(I-F)^{-1}G+J)\mathcal{F}$ in Theorem~\ref{thm:circper} is not necessarily diagonal \cite{ChenQiu1997}.
\end{rem}
\begin{pf}[Theorem~\ref{thm:circper}]
Substituting $a_k$ from \eqref{eq:a_k}, letting $\omega_m = \frac{2 \pi m}{N}$ and using
\begin{math}
\sum_{k=1}^N M^k
= (I - M)^{-1} (M - M^{N+1}) 
= (I - M)^{-1} (I - M^N) M,
\end{math}
leads to
\begin{math}
\lambda_m
= \sum_{k=0}^{N-1} a_k e^{-j \omega_m k} 
= D + \sum_{k=0}^{N-1} C A^{N-k-1} (I - A^N)^{-1} B e^{-j \omega_m k} 
\hspace{-3.7mm}\overset{\tiny \tilde{k} := N-k}{=} D + \sum_{\tilde{k}=1}^{N} C A^{\tilde{k}-1} e^{j \omega_m \tilde{k}} (I - A^N)^{-1} B 
= D + C \left( \sum_{\tilde{k}=1}^{N} (A e^{j \omega_m})^{\tilde{k}} \right) (I - A^N)^{-1} A^{-1} B 
= D \; + 
\quad C (I - A e^{j \omega_m})^{-1} (I - A^N) A e^{j \omega_m} (I - A^N)^{-1} A^{-1} B 
= D + C (e^{-j \omega_m} I - A)^{-1} B
= P(e^{-j \omega_m}),
\end{math}
which yields the desired result. \qed
\end{pf}

If $P \in \mathcal{R}\mathcal{H}_\infty$, then from Theorem~\ref{thm:circper} it follows that $\lim_{N \rightarrow \infty}\bar \sigma(H(I-F)^{-1}G+J) \rightarrow \|P\|_\infty$. Compared to \cite[Theorem 3]{RojasOomHjaWah2012}, the result in Theorem~\ref{thm:circper} only involves frequency discretization, and hence transient effects do not play a role. This improves norm estimation properties for finite $N$, as is confirmed in \Sec~\ref{sec:example}.

\section{The need for time reversal: Reversed-circulant matrices}

Motivated by Theorem~\ref{thm:circper}, the central idea in this paper is to apply the power iteration to the matrix $H(I-F)^{-1}G+J$.  In the periodically operated case, there is a direct connection between the eigenvalues and the frequency response, and hence properties such as the $\Hinf$ norm. However, the input $u \in \mathbb{R}^N$ must be real-valued, while the eigenvalues of $H(I-F)^{-1}G+J$ are complex by Theorem~\ref{thm:circper}. These complex eigenvalues introduce a rotation and hence convergence cannot be expected if the power iteration is applied directly to $H(I-F)^{-1}G+J$, see \cite[\Sec~7.12]{Wilkinson1965}. 

The following theorem is key to the subsequent steps.
\begin{thm}\label{thm:eqaftertimreversal}
Let $\mathrm{circ}(a_0, a_1, \ldots a_{N-1})$ be a given circulant matrix with spectrum $\lambda_m, m = 0, 1, \ldots, N-1$ in \eqref{eq:dft}. Then, 
\begin{align}
& R = \tijdreverse\mathrm{circ}(a_0, a_1, \ldots a_{N-1}) =\\
& \quad 
\begin{bmatrix}
a_1 & a_2 & a_3 & \cdots & a_{N-1} & a_{0} \\
a_{2} & a_3 & a_4 & \cdots & a_{0} & a_{1}\\  
a_{3} & a_{4} & a_5 & \ddots & a_{1} & a_{2}\\
\vdots & \vdots & \ddots & \ddots & \ddots & \vdots \\
a_{N-1} & a_0 & a_1 & \ddots & a_{N-3} & a_{N-2} \\
a_0 & a_1 & a_2 & \cdots & a_{N-2} & a_{N-1}
\end{bmatrix}
\end{align}
is a reversed-circulant matrix with eigenvalues given by 
\begin{align}
\lambda^{\tijdreverse}_0  &= \lambda_0\\
\lambda^{\tijdreverse}_{m} &= -\lambda^{\tijdreverse}_{N-m}  = \mathrm{abs}(\lambda_m), \quad 0 < m < N/2 \\
\lambda^{\tijdreverse}_{N/2} &= -\lambda_{N/2}.
\end{align}
\end{thm}
\begin{pf}
Let $C = \mathrm{circ}(a_0, a_1, \ldots a_{N-1})$, and note that $R = \sum_{k=0}^{N-1} a_k \tijdreverse Y^k$, where $Y \in \mathbb{R}^{N \times N}$ satisfies $Y [x_1, x_2, \dots, x_N]^T = [x_2, \dots, x_N, x_1]^T$ for $x_1, \dots, x_N \in \mathbb{C}$.
To establish the eigenstructure of $R$, let $w_m := [1, W_N^m, \dots, W_N^{(N-1) m}]^T \in \mathbb{R}^N$ for $m = 0, 1, \dots, N-1$. As noted in the proof of Theorem~\ref{thm:circper}, $C$ has $\lambda_0, \dots, \lambda_{N-1}$ as eigenvalues, and $w_0, \dots, w_{N-1}$ as the corresponding eigenvectors, which are thus linearly independent.

For $m=0$, we have that $w_0 = [1,\dots,1]^T$, so $R w_0 = \tijdreverse C w_0 = \lambda_0 \tijdreverse w_0 = \lambda_0 w_0$, hence $\lambda_0$ is an eigenvalue of $R$ with corresponding eigenvector $w_0$.

If $N$ is even, then $w_{N/2} = [1,-1, \dots, 1, -1]^T$, so $\tijdreverse w_{N/2} = - w_{N/2}$. Therefore, $R w_{N/2} = \tijdreverse C w_{N/2} = \lambda_{N/2} \tijdreverse w_{N/2} = -\lambda_{N/2} w_{N/2}$. Hence, $-\lambda_{N/2}$ is an eigenvalue of $R$ with corresponding eigenvector $w_{N/2}$.

For $0 < m < N/2$, first assume that $\lambda_m \neq 0$, and let $v_{\pm m} := |\lambda_m| w_m \pm \lambda_m W_N^{-m} w_{N-m}$. The vectors $v_m$ and $v_{-m}$ are linearly independent (and also from the remaining vectors $w_m$, $m \notin \{ m, N - m \}$, as they are built only from $w_m$ and $w_{N-m}$), since $[v_m \; v_{-m}] = [w_m \; w_{N-m}]\, A$, where $A := \begin{bmatrix} |\lambda_m| & \lambda_m W_N^{-m} \\ |\lambda_m| & -\lambda_m W_N^{-m} \end{bmatrix}$ has determinant $- 2 \lambda_m |\lambda_m| W_N^{-m} \neq 0$, since $\lambda_m \neq 0$.

Notice that $\lambda_{N-m} = \lambda_{-m} = \overline{\lambda_m}$, and that $Y^k w_m = W_N^{m k} w_m$. Then,
\begin{math}
R v_{\pm m} 
= |\lambda_m| R w_m \pm \lambda_m W_N^{-m} R w_{N-m} 
= |\lambda_m| \sum_{k=0}^{N-1} a_k \tijdreverse Y^k w_m \pm \lambda_m W_N^{-m} \sum_{k=0}^{N-1} a_k \tijdreverse Y^k w_{N-m} \\ 
= |\lambda_m| \sum_{k=0}^{N-1} a_k \tijdreverse W_N^{m k} w_m \pm \lambda_m W_N^{-m} \sum_{k=0}^{N-1} a_k \tijdreverse W_N^{-m k} \\ w_{N-m} 
= |\lambda_m| \lambda_m \tijdreverse w_m \pm |\lambda_m|^2 W_N^{-m} \tijdreverse w_{N-m}.
\end{math}
For all $l = 1, \dots, N$, $[\tijdreverse w_m]_l = [w_m]_{N+1-l} = W_N^{(N-l)m}$, so
\begin{math}
[v_{\pm m}]_l
= [|\lambda_m| w_m \pm \lambda_m W_N^{-m} w_{N-m}]_l 
= |\lambda_m| W_N^{(l-1) m} \pm \lambda_m W_N^{-m} W_N^{(l-1)(N-m)} 
= |\lambda_m| W_N^{(l-1) m} \pm \lambda_m W_N^{-l m},
\end{math}
while
\begin{math}
[R v_{\pm m}]_l
= [|\lambda_m| \lambda_m \tijdreverse w_m \pm |\lambda_m|^2 W_N^{-m} \tijdreverse w_{N-m}]_l 
= |\lambda_m| \lambda_m W_N^{(N-l)m} \pm |\lambda_m|^2 W_N^{-m} W_N^{(N-l)(N-m)} \\ 
= |\lambda_m| \lambda_m W_N^{-lm} \pm |\lambda_m|^2 W_N^{(l-1) m} 
= \pm |\lambda_m| (|\lambda_m| W_N^{(l-1) m} \pm \lambda_m W_N^{-l m})
\end{math}
hence $R v_{\pm m} = \pm |\lambda_m| v_{\pm m}$, so $\pm |\lambda_m|$ is an eigenvalue of $R$ with eigenvector $v_{\pm m}$.

Consider all the vectors $w_m$ for which $\lambda_m = 0$. These are eigenvectors of $C$ corresponding to the eigenvalue $0$, \emph{i.e.}, $C w_m = 0$. Then, $R w_m = \tijdreverse C w_m = 0$, so $w_m$ is also an eigenvector of $R$ corresponding to the eigenvalue $\lambda_m = 0$. Together with the previous eigenvectors, they span $\mathbb{R}^N$, which thus exhausts the search for eigenvectors of $R$. This concludes the proof. \qed
\end{pf}
Hence, Theorem~\ref{thm:eqaftertimreversal} reveals that time-reversal essentially enforces the matrix $\tijdreverse(H(I-F)^{-1}G+J)$ to be symmetric, implying real-valued eigenvalues that have the same magnitude as those of $H(I-F)^{-1}G+J$. 

\section{Algorithm}\label{sec:alg}

In this section, the power iteration is implemented. Additional steps are taken to actually implement the algorithm on the matrix $\tijdreverse(H(I-F)^{-1}G+J)$ to enforce the periodic response~\eqref{eq:periodicresponse} and bounded inputs. To this end, let
\begin{equation}\label{eq:controller}
u_{j+1} = \alpha_k \mu_j \left(\tijdreverse y_j  + \lambdashift I_N \right)+ (1 - \alpha_k) u_j,
\end{equation}
where $\lambdashift \in \mathbb{R}\backslash 0$ is discussed in Lemma~\ref{lemma:PIM}, and
\begin{equation}\label{eq:normalization}
\mu_j = \frac{\|y_j\|_2}{\sqrt{N}}
\end{equation}
is a normalization constant to normalize the input power to one, see \cite{RojasOomHjaWah2012}. In addition, 
\begin{equation}\label{eq:update}
\alpha_j = 
\begin{cases}
1 & \text{if } j = n^{\mathrm{update}}q, q \in \mathbb{Z}_{>0}\\
0 & \text{otherwise}.
\end{cases}
\end{equation}
Here, $\alpha_j$ enforces that the input is only updated every $n^{\mathrm{update}}$ iterations to ensure that a steady-state response is obtained. Note that this additional transient time is at least of the same order of magnitude as approaches based on resetting, including \cite{RojasOomHjaWah2012}, since it takes a similar transient time to bring a system to rest. Practically, $n^{\mathrm{update}}$ can be set such that the system is sufficiently in steady-state, which is similar to, e.g., periodic excitation in frequency domain system identification or the iterative inversion-based control, see \cite{RozarioOom2023} and references therein for details. 

\begin{lemma}\label{lemma:PIM}
Let $u_1$ be a random vector. Then, the vectors $u_j$ generated by iteration \eqref{eq:controller} converge to the eigenvector corresponding to the largest eigenvalue $\lambda_m+\lambdashift$, where $\lambda_m$ is defined in \eqref{eq:dft}. 
\end{lemma} 
A proof of Lemma~\ref{lemma:PIM} follows from \cite[\Sec\ 9.3]{Wilkinson1965}, where the randomness of the vector $u_1$ is used to ensure it has a component in the direction of the eigenvector corresponding to $\lambda_m+\lambdashift$. The shift $\lambdashift$ is essential to let the algorithm converge to the largest positive eigenvalue in Theorem~\ref{thm:eqaftertimreversal}, since there is always a negative one with the same magnitude. Indeed, a negative eigenvalue would violate the input in \eqref{eq:convinput}, since the input would start switching signs. A suitable order of magnitude for the choice regarding $\lambdashift$ is given by $\|P\|_\infty$, which renders all eigenvalues $\lambda_m$ positive, enhancing convergence, see \cite[Remark 9]{RojasOomHjaWah2012}. If there are multiple dominant eigenvectors $\lambda_m+\lambdashift$, then the eigenvector iteration will converge to a subspace that exhibits similar norm estimation properties \cite[\Sec\ 9.3]{Wilkinson1965}. Finally, it is noted that standard results apply regarding conditioning of the matrix \cite[\Sec\ 9.3]{Wilkinson1965}, yet it is emphasized that the presence of noise may be dominant in case the iteration \eqref{eq:controller} is applied to a physical system, see \cite{RojasOomHjaWah2012} for a detailed analysis.

Given the iteration~\eqref{eq:controller}, an estimator for the $\Hinf$ norm is directly obtained from \eqref{eq:Ri2_fin_time}, see also \cite{RojasOomHjaWah2012}, as
\begin{math}
\beta_\iter = \tfrac{1}{N}{u_\iter^T y_\iter}.
\end{math}

\section{Example}\label{sec:example}

\begin{figure}
\centering
\input{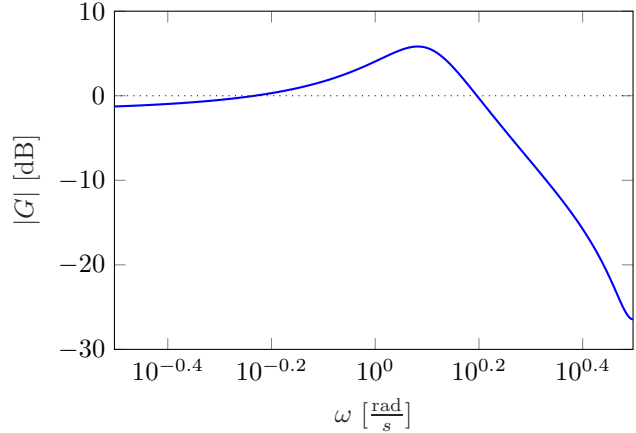}
\caption{Bode magnitude plot of $P$ in \eqref{eq:exsyst}.}
\label{fig:bodeplot}
\end{figure}

Consider the system 
\begin{equation}\label{eq:exsyst}
P(z) = z^{-50}\frac{5 z^{-1} + 4 z^{-2}}{10 - 5 z^{-1} + 6 z^{-2}}, 
\end{equation}
see \Fig~\ref{fig:bodeplot} for a Bode magnitude plot. 

\begin{figure}
\centering
%
%
\begin{tikzpicture}

\begin{axis}[%
width=0.951\figurewidth,
height=\figureheight,
at={(0\figurewidth,0\figureheight)},
scale only axis,
xmode=log,
xmin=1,
xmax=512,
xtick={1,2,4,8,16,32,64,128,256,512},
xticklabels={{$\text{2}^\text{1}$},{$\text{2}^{\text{1}}$},{$\text{2}^{\text{2}}$},{$\text{2}^{\text{3}}$},{$\text{2}^{\text{4}}$},{$\text{2}^{\text{5}}$},{$\text{2}^{\text{6}}$},{$\text{2}^{\text{7}}$},{$\text{2}^{\text{8}}$},{$\text{2}^{\text{9}}$}},
xminorticks=true,
xlabel style={font=\color{white!15!black}},
xlabel={$N$},
ymin=0.5,
ymax=2,
ytick={0.5,   1, 1.5,   2},
ylabel style={font=\color{white!15!black}},
ylabel={Norm estimate},
axis background/.style={fill=white}
]
\addplot [color=blue, line width=0.8pt, forget plot]
  table[row sep=crcr]{%
1	0.818181818181818\\
2	0.818181818181818\\
4	1\\
8	1.20314440034967\\
16	1.94280657856413\\
64	1.94280657856413\\
128	1.94931220339798\\
256	1.95440518251759\\
512	1.95440518251759\\
};
\addplot [color=red, line width=0.8pt, forget plot]
  table[row sep=crcr]{%
51.6979584750314	0.35\\
52	0.5\\
53	0.921343022093828\\
54	1.04291660013137\\
55	1.06327164423811\\
56	1.24763122559185\\
57	1.27852100892525\\
58	1.36852323860787\\
59	1.4394219463788\\
60	1.4444461811811\\
61	1.53633995360825\\
62	1.55823317468255\\
63	1.59358114537293\\
64	1.63629072055846\\
65	1.63832402168906\\
66	1.68390099325265\\
67	1.70050023461134\\
68	1.7132706669477\\
69	1.74121948296756\\
70	1.74499861505543\\
71	1.76656698210475\\
72	1.77954182714269\\
73	1.78304189888315\\
74	1.80234382799541\\
75	1.80671283305469\\
77	1.82714128866176\\
78	1.82727003381276\\
79	1.84066036060794\\
81	1.84946061263328\\
82	1.85773978288034\\
83.0000000000001	1.85872116901827\\
84	1.86613413934895\\
85	1.87031033073839\\
86	1.87179601494078\\
87	1.87847025409374\\
88	1.87993829393624\\
90	1.88772066277069\\
91.0000000000001	1.88774893015384\\
92	1.89312936901298\\
94	1.89676738529744\\
95.0000000000001	1.90019468750398\\
96	1.90057298894234\\
97	1.90386970932164\\
98.0000000000001	1.9056816718431\\
99	1.90640703734301\\
100	1.90942175266042\\
101	1.91007597233774\\
103	1.91378070417764\\
104	1.9138253414945\\
105	1.91643490076502\\
107	1.9182578719128\\
108	1.91998299647304\\
109	1.92016655016805\\
110	1.92189208398019\\
111	1.92283235056049\\
112	1.92323089260956\\
113	1.9248337863839\\
114	1.92518006586353\\
116	1.92721290924696\\
117	1.92724590373147\\
118	1.92869897153332\\
120	1.92973430344374\\
121	1.93071975565021\\
122	1.93082298232121\\
123	1.93183087024743\\
129	1.93500408267264\\
130	1.93502620244395\\
131	1.93591569588904\\
135	1.93723637352424\\
137	1.93822123784717\\
138	1.93837418261563\\
139	1.93898290766934\\
141	1.93952506922804\\
145	1.94070089298159\\
171	1.94590166958368\\
210	1.94959962244207\\
279	1.95220521818791\\
358	1.95333774720097\\
487	1.95405205514834\\
512	1.95412589085292\\
};
\addplot [color=black, line width=0.8pt, forget plot]
  table[row sep=crcr]{%
1	1.95477063458545\\
512	1.95477063458545\\
};
\end{axis}
\end{tikzpicture}%
\caption{Norm estimate comparison $\lambdamax(\mathcal{F}^{*}(H(I-F)^{-1}G+J)\mathcal{F})$ (\Colorone), $\lambdamax(\tijdreverse J)$ (\Colortwo), and $\|P\|_\infty$ (\Colorthree). The estimate $\lambdamax(\mathcal{F}^{*}(H(I-F)^{-1}G+J)\mathcal{F})$ converges much faster to the true system norm $\|P\|_\infty$, hence the reset-free approach leads to a better approximation.}
\label{fig:normestimate}
\end{figure}
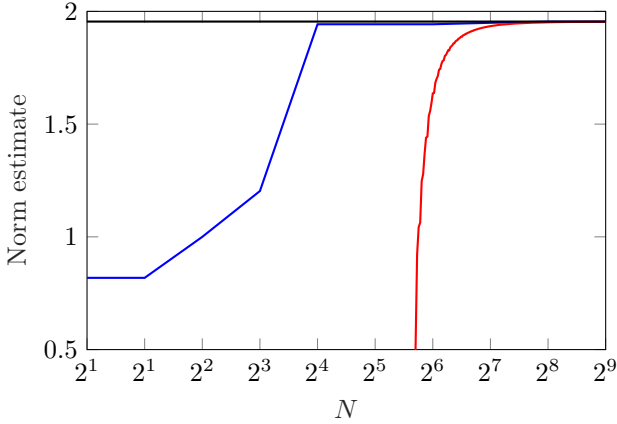

First, the result of Theorem~\ref{thm:circper} is compared to the traditional reset-based appproach~\eqref{eq:Ri2_fin_time3}. In particular, the spectral radius $\lambdamax(\mathcal{F}^{*}(H(I-F)^{-1}G+J)\mathcal{F})$ is computed for different values of $N$. In addition, the traditional result based on resetting is computed, i.e., $\lambdamax(\tijdreverse J)$. The values are computed for a range of values of $N$. For $\lambdamax(\mathcal{F}^{*}(H(I-F)^{-1}G+J)\mathcal{F})$, $N$ is increased by doubling its value, ensuring that the frequency points are retained. Furthermore, these approximations for finite $N$ are compared to their true value $\|P\|_\infty$. The results are presented in \Fig~\ref{fig:normestimate}. Clearly, the value of $\lambdamax(\mathcal{F}^{*}(H(I-F)^{-1}G+J)\mathcal{F})$ converges much faster to the true value $\|P\|_\infty$ as expected. In addition, the result $\lambdamax(\tijdreverse J)$ does produce a norm estimate of zero for $N<50$ due to the strict delays in \eqref{eq:exsyst}. The reset-free estimate $\lambdamax(\mathcal{F}^{*}(H(I-F)^{-1}G+J)\mathcal{F})$ clearly does not suffer from these delays. 

Next, the approach of \Sec~\ref{sec:alg} is implemented. Here, $N = 50$ and $n^{\mathrm{update}} = 10$. The resulting input signals are depicted in \Fig~\ref{fig:inputs}. The input already after one update contains a dominant frequency component, which further becomes sinusoidal after several iterations. This frequency excites the resonance peak in \Fig~\ref{fig:bodeplot}, without the transient envelope effects that typically arise in the reset-based approach of \cite{RojasOomHjaWah2012}, see \cite[\Fig~11]{OomenMaaRojHja2014}. This is also evidenced in \Fig~\ref{fig:normestimateiterative}, where the estimated norm $\beta_j$ rapidly converges to the underlying value $\lambdamax(\mathcal{F}^{*}(H(I-F)^{-1}G+J)\mathcal{F})$. Also, the update of the input each $n^{\mathrm{update}} = 10$ iterations can be distinguished, leading to a transient in the iteration domain after each update.

Summarizing, the approach of \Sec~\ref{sec:alg} effectively estimates the norm directly from data through an iterative procedure. To further compare this procedure with pre-existing approaches, including the reset-based approach in \cite{RojasOomHjaWah2012}, note that $J$ in \eqref{eq:liftedsystem} is equal to $J=0$ for $N=50$ due to the strict delay in \eqref{eq:exsyst}. Consequently, a reset-based approach does not lead to a useful estimate of the gain. 

The superior estimation results of the reset-free algorithm are promising for robust control and related fields. In particular, robust control typically require an accurate model error bound in terms of the $\mathcal{H}_\infty$ norm. The reset-free approach can be directly applied to the model error, see \cite{OomenMaaRojHja2014}, leading to very accurate results with a reduced experimental burden to reset the system. Similarly, the algorithm enables online fault detection and isolation through suitable experiments. Finally, the algorithm can be seen as an online optimal experiment design, providing more accurate models in system identification as well as related data-driven estimation algorithms.  

\setlength\figureheight{9cm}
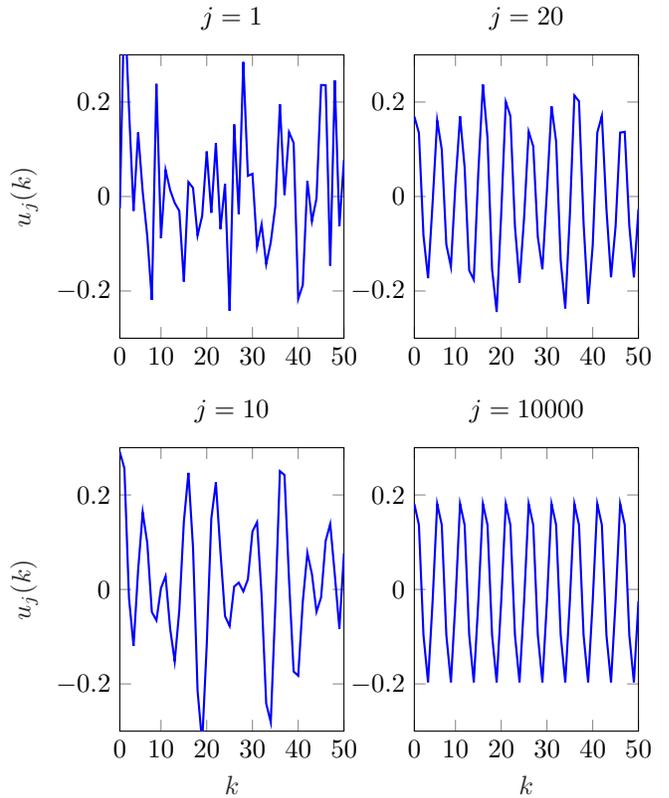
\begin{figure}
\centering
%
%
\begin{tikzpicture}

\begin{axis}[%
width=0.411\figurewidth,
height=0.419\figureheight,
at={(0\figurewidth,0.581\figureheight)},
scale only axis,
unbounded coords=jump,
xmin=1,
xmax=50,
xtick={1,10,20,30,40,50},
xticklabels={{0},{10},{20},{30},{40},{50}},
ymin=-0.3,
ymax=0.3,
ylabel style={font=\color{white!15!black}},
ylabel={$u_j(k)$},
axis background/.style={fill=white},
title style={font=\bfseries},
title={$j=1$}
]
\addplot [color=blue, line width=0.8pt, forget plot]
  table[row sep=crcr]{%
1	-0.0248918966113933\\
1.86123365395217	0.359999999999999\\
nan	nan\\
2.26370693792216	0.359999999999999\\
3	0.186846577300493\\
4	-0.0303673061294063\\
5	0.135657156372041\\
6	0.0101227733649267\\
7	-0.0844381690795188\\
8	-0.218689088564844\\
9	0.238398107508175\\
10	-0.0882753301697932\\
11	0.05613760283277\\
12	0.0137424094315435\\
13	-0.0134965230233348\\
14	-0.0304878930672956\\
15	-0.180050323347004\\
16	0.0305481062897641\\
17	0.0183242274164428\\
18	-0.0826101666442511\\
19	-0.0409435018096005\\
20	0.095347846788485\\
21	-0.0344023547431149\\
22	0.11319479243695\\
23	-0.0688397428765057\\
24	0.0271200886838727\\
25	-0.241954011842331\\
26	0.152428110669923\\
27	-0.0370431390489614\\
28	0.284822836037762\\
29	0.0437482447929156\\
30	0.0480373057513717\\
31	-0.106323041790347\\
32	-0.0590093730489514\\
33	-0.143215633253661\\
34	-0.0969062444852185\\
35	-0.0204012485813507\\
36	0.195315100487448\\
37	0.00329229103044071\\
38	0.136023480280627\\
39	0.11352850665002\\
40	-0.216769583843643\\
41	-0.187431020018245\\
42	0.0330923653089812\\
43	-0.052140944789663\\
44	-0.00388262258395855\\
45	0.236012320012449\\
46	0.235826724781475\\
47	-0.146357476173613\\
48	0.245535086333199\\
49	-0.0629221202582713\\
50	0.0779824865697165\\
};
\end{axis}

\begin{axis}[%
width=0.411\figurewidth,
height=0.419\figureheight,
at={(0\figurewidth,0\figureheight)},
scale only axis,
xmin=1,
xmax=50,
xtick={1,10,20,30,40,50},
xticklabels={{0},{10},{20},{30},{40},{50}},
xlabel style={font=\color{white!15!black}},
xlabel={$k$},
ymin=-0.3,
ymax=0.3,
ylabel style={font=\color{white!15!black}},
ylabel={$u_j(k)$},
axis background/.style={fill=white},
title style={font=\bfseries},
title={$j=10$}
]
\addplot [color=blue, line width=0.8pt, forget plot]
  table[row sep=crcr]{%
1	0.290839295494649\\
2	0.256298612095847\\
3	-0.0179470038780138\\
4	-0.119675089419793\\
5	0.040896494241828\\
6	0.163250882440693\\
7	0.0983738546464252\\
8	-0.0477716594297917\\
9	-0.0654930698157727\\
10	0.00272977803166441\\
11	0.0279605052404506\\
12	-0.0829267959558138\\
13	-0.152092637453023\\
14	-0.0443369119009063\\
15	0.14376031037483\\
16	0.246758564342791\\
17	0.0892771539482879\\
18	-0.214760438804028\\
19	-0.324256299573655\\
20	-0.117771056564102\\
21	0.15006982209853\\
22	0.226506053139154\\
23	0.0791214338438735\\
24	-0.0567015577898573\\
25	-0.0774072958558847\\
26	0.00609537021178852\\
27	0.0143051885070804\\
28	-0.00421094137401923\\
29	0.0215050430552495\\
30	0.122823295171827\\
31	0.142212444166496\\
32	-0.0134643418846494\\
33	-0.241487727130377\\
34	-0.28119193192596\\
35	-0.0339737354319567\\
36	0.250426553721468\\
37	0.242259412999601\\
38	0.0188542847006161\\
39	-0.173848033648731\\
40	-0.182672501352961\\
41	-0.0258315690117996\\
42	0.0771280268244965\\
43	0.032428346072912\\
44	-0.0473427527620558\\
45	-0.0169325677650889\\
46	0.102993491091105\\
47	0.13856813535152\\
48	0.0309184189133944\\
49	-0.0830642988592416\\
50	0.0763203859075219\\
};
\end{axis}

\begin{axis}[%
width=0.411\figurewidth,
height=0.419\figureheight,
at={(0.54\figurewidth,0.581\figureheight)},
scale only axis,
xmin=1,
xmax=50,
xtick={1,10,20,30,40,50},
xticklabels={{0},{10},{20},{30},{40},{50}},
ymin=-0.3,
ymax=0.3,
axis background/.style={fill=white},
title style={font=\bfseries},
title={$j=20$}
]
\addplot [color=blue, line width=0.8pt, forget plot]
  table[row sep=crcr]{%
1	0.169545754330592\\
2	0.134779172021837\\
3	-0.0825686524940394\\
4	-0.172907888617473\\
5	-0.0129126486214375\\
6	0.161054499640038\\
7	0.0998026950755175\\
8	-0.101891775972234\\
9	-0.148592011748072\\
10	0.0276827247059686\\
11	0.169506124117525\\
12	0.061645753139338\\
13	-0.15662417398093\\
14	-0.176051168153833\\
15	0.0530572011318142\\
16	0.237147907423854\\
17	0.124208527132609\\
18	-0.15364434348011\\
19	-0.244453041663284\\
20	-0.0360562500771309\\
21	0.200248677820888\\
22	0.170353836385004\\
23	-0.0631274062041598\\
24	-0.182146995779775\\
25	-0.0450257761924604\\
26	0.13741450054264\\
27	0.106376025721055\\
28	-0.0856106029273178\\
29	-0.154080906274586\\
30	0.0138597742588047\\
31	0.190940128111286\\
32	0.116434649710008\\
33	-0.13568162856825\\
34	-0.236754192921062\\
35	-0.0395191988261843\\
36	0.213877819064912\\
37	0.201721831212161\\
38	-0.0561915298031579\\
39	-0.227276217503523\\
40	-0.103191924574425\\
41	0.135331572150122\\
42	0.170581083782174\\
43	-0.0266552256259871\\
44	-0.170935695563543\\
45	-0.0632620228924523\\
46	0.135407805238586\\
47	0.136927325659002\\
48	-0.0615646610563516\\
49	-0.171054067562892\\
50	-0.0263248460142691\\
};
\end{axis}

\begin{axis}[%
width=0.411\figurewidth,
height=0.419\figureheight,
at={(0.54\figurewidth,0\figureheight)},
scale only axis,
xmin=1,
xmax=50,
xtick={1,10,20,30,40,50},
xticklabels={{0},{10},{20},{30},{40},{50}},
xlabel style={font=\color{white!15!black}},
xlabel={$k$},
ymin=-0.3,
ymax=0.3,
axis background/.style={fill=white},
title style={font=\bfseries},
title={$j=10000$}
]
\addplot [color=blue, line width=0.8pt, forget plot]
  table[row sep=crcr]{%
1	0.181096016218106\\
2	0.136685286455396\\
3	-0.0966198634266533\\
4	-0.196399646041442\\
5	-0.0247617932054069\\
6	0.181096016218106\\
7	0.136685286455396\\
8	-0.0966198634266533\\
9	-0.196399646041442\\
10	-0.0247617932054069\\
11	0.181096016218106\\
12	0.136685286455396\\
13	-0.0966198634266533\\
14	-0.196399646041442\\
15	-0.0247617932054069\\
16	0.181096016218106\\
17	0.136685286455396\\
18	-0.0966198634266533\\
19	-0.196399646041442\\
20	-0.0247617932054069\\
21	0.181096016218106\\
22	0.136685286455396\\
23	-0.0966198634266533\\
24	-0.196399646041442\\
25	-0.0247617932054069\\
26	0.181096016218106\\
27	0.136685286455396\\
28	-0.0966198634266533\\
29	-0.196399646041442\\
30	-0.0247617932054069\\
31	0.181096016218106\\
32	0.136685286455396\\
33	-0.0966198634266533\\
34	-0.196399646041442\\
35	-0.0247617932054069\\
36	0.181096016218106\\
37	0.136685286455396\\
38	-0.0966198634266533\\
39	-0.196399646041442\\
40	-0.0247617932054069\\
41	0.181096016218106\\
42	0.136685286455396\\
43	-0.0966198634266533\\
44	-0.196399646041442\\
45	-0.0247617932054069\\
46	0.181096016218106\\
47	0.136685286455396\\
48	-0.0966198634266533\\
49	-0.196399646041442\\
50	-0.0247617932054069\\
};
\end{axis}
\end{tikzpicture}%
\caption{Input signals. The initial input signal $u_1(k)$ (\Colorone) is random noise. Already after the first update, the input $u_{10}(k)$ (\Colortwo) contains a dominant frequency component. After convergence, the input, $u_{10000}(k)$ (\Colorthree) clearly contains a dominant frequency component that corresponds to the peak frequency in \Fig~\ref{fig:bodeplot}.}
\label{fig:inputs}
\end{figure} 
\setlength\figureheight{4.5cm} 

\begin{figure}
\centering
\input{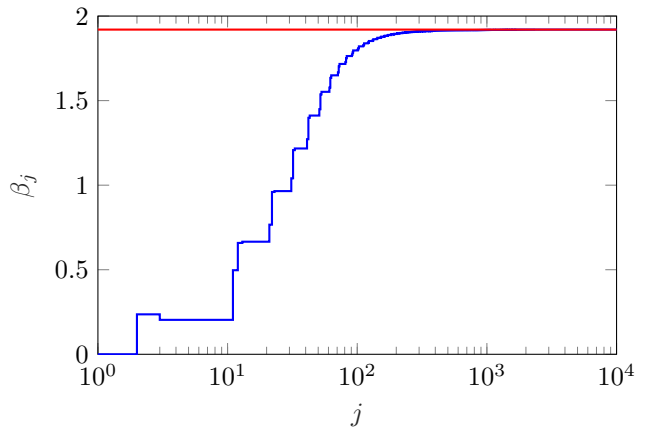}
\caption{Estimated norm $\beta_j$ using the approach of \Sec~\ref{sec:alg} (\Colorone) and $\lambdamax(\mathcal{F}^{*}(H(I-F)^{-1}G+J)\mathcal{F})$ for $N=50$ (\Colortwo). Here, $\beta_j$ converges rapidly to the underlying value  $\lambdamax(\mathcal{F}^{*}(H(I-F)^{-1}G+J)\mathcal{F})$.}
\label{fig:normestimateiterative}
\end{figure}

\section{Conclusion}
In this paper, the standard assumption regarding resetting the system at each iteration in iterative data-driven estimation algorithms is alleviated. This has both theoretical advantages in terms of more accurate estimates and practical advantages in terms of experimental effort. Its theoretical advantages are illustrated on a data-driven $\Hinf$ norm estimation problem. The practical advantages compared to earlier data-driven norm estimation algorithms are reduced experimental time and reduced user intervention, since it does not require cumbersome resetting and re-initialisation procedures of the experimental setup. In this respect, the reset-free method has major practical advantages compared to earlier algorithms and no foreseen disadvantages. The extension to related data-driven estimation algorithms is immediate, and may enhance theoretical and practical advantages of a large class of these approaches. \\

For future research, it is of interest to explore whether it is possible to reduce $n^{\mathrm{update}}$ in \eqref{eq:update}. Currently, convergence is only guaranteed if $n^{\mathrm{update}}$ is selected such that the input reaches a steady-state. Numerous experiments have shown that the algorithm works practically if the input in iteration $j$ has limited leakage to iteration $j+1$, in which case $n^{\mathrm{update}}=1$ often leads to convergence. A formal convergence proof is complicated by the normalization~\eqref{eq:normalization}, see also \cite{RojasOomHjaWah2012}. Furthermore, bounds on the error $\beta_j - \|P\|_\infty$ depending on $N$ are of interest for selecting $N$ based on prior knowledge. 



\end{document}